\documentclass{amsproc}
\usepackage{amssymb}
\usepackage{amsfonts}

\theoremstyle{plain}

\newtheorem{corollary}{Corollary}

\newtheorem{example}{Example}

\newtheorem{lemma}{Lemma}

\newtheorem{proposition}{Proposition}
\newtheorem{remark}{Remark}

\newtheorem{theorem}{Theorem}
\numberwithin{equation}{section}

\input{tcilatex}

\begin{document}
\title[Iterated symmetric bimodal maps]{Irreducible complexity of iterated
symmetric bimodal maps}
\author{J. P. Lampreia }
\address[J. P. Lampreia ]{Department of Mathematics, F.C.T. \\
Universidade Nova de Lisboa, Lisboa, Portugal}
\email[J. P. Lampreia]{jpl@fct.unl.pt}
\author{R. Severino }
\address[R. Severino ]{Department of Mathematics\\
Universidade do Minho Braga, Portugal}
\email[R. Severino]{ricardo@math.uminho.pt}
\author{J. Sousa Ramos}
\address[J. Sousa Ramos]{Department of Mathematics, I.S.T.\\
Universidade T\'{e}cnica de Lisboa, Lisboa, Portugal}
\email[J. Sousa Ramos]{sramos@math.ist.utl.pt}
\urladdr{http://www.math.ist.utl.pt/\~{ }sramos}
\date{February, 5, 2004}
\subjclass[2000]{Primary 37E05, 37B10; Secondary 37E20, 37B99}
\thanks{The authors gratefully acknowledge the financial support from the
project FCT-POCTI-FEDER}

\begin{abstract}
We introduce a tree structure for the iterates of symmetric bimodal maps and
identify a subset which we prove to be isomorphic to the family of unimodal
maps. This subset is used as a second factor for a $\ast $-product that we
define in the space of bimodal kneading sequences. Finally, we give some
properties for this product and study the $*$-product induced on the
associated Markov shifts.
\end{abstract}

\maketitle


\section{Introduction and preliminary definitions}

The concept of irreducible complexity of a biological system was introduced
by Behe, \cite{Behe}, in 1996. His point of view is that an organim
consisting of a finite, possibly very large, number of independent
components, coupled together in some way, exhibits irreducible complexity
if, by removing any of its component, the reduced system no longer functions
meaningfully. Using the language of non-linear dynamics and chaos theory,
Boyarsky and G\'ora, \cite{BoyarskyGora}, reinterpreted Behe's definition
from a Markov transition matrix perspective by saying that a system is
irreducibly complex if the associated transition matrix is primitive but no
principal submatrix is primitive.

It is our conviction that the concept of reducible complexity of a dynamical
system can also be interpreted in terms of a factorization: within Milnor
and Thurston's kneading theory framework, and the topological classification
obtained from it, Derrida, Gervois, and Pomeau, \cite{De-Ge-Po 78},
introduced a $\star$-product between unimodal kneading sequences for which
it was possible to prove that the topological entropy, a measure of
complexity, of a factorizable system is equal to the topological entropy of
one of the factors. Despite of a larger number of its components, the
complexity of the system remains the same whenever its irreducible
component, a factor of the product, does not change.

Some years latter, Lampreia, Rica da Silva, and Sousa Ramos, \cite{La-RS-SR
88}, introduced a Markov transition matrix formalism associated with the
kneading theory and a product between unimodal matrices corresponding to the
Derrida, Gervois, and Pomeau $\star$-product. Then, they proved that
irreducible unimodal kneading sequences corresponds to primitive Markov
transition matrices.

With this work we would like to introduce the generalization, for bimodal
symmetric maps of the interval, of the $\star$-product and the corresponding
product between transition matrices.

Consider a two-parameter family $f_{a,b}$ of maps, from the closed interval $%
I=[c_0,c_3]$ into itself, with two critical points, usually called a bimodal
family of maps of the interval, see \cite{Co-Ec 80}, \cite{Mi-Th 88}, \cite%
{La-SR 97}. Once fixed the parameters $(a,b)$, the map $f_{a,b}$ is
piecewise monotone and hence $I$ can be subdivided in the following three
subintervals: $L=[c_{0},c_{1}]$, $M=[c_{1},c_{2}]$ and $R=[c_{2},c_{3}]$,
where $c_{i}$ are the critical points or the extremal points, in such a way
that the restriction of $f$ to each interval is strictly monotone. We will
choose the family of maps such that the restrictions $f_{a,b|L}$ and $%
f_{a,b|R}$ are increasing and the restriction $f_{a,b|M}$ is decreasing.

For each value $(a,b)$ we define the orbits of the critical points by: 
\begin{equation*}
O(c_{i})=\{x_{j}:x_{j}=f^{j}(c_{i}),~j\in \mathbb{N}\}
\end{equation*}%
with $i=1,2$.

With the aim of studying the topological properties of these orbits we
associate to each orbit $O(c_{i})$ a sequence of symbols $S=S_{1}S_{2}\dots
S_{j}\dots $ where $S_{j}=L$ if $f_{a,b}^{j}(c_{i})<c_{1}$, $S_{j}=A$ if $%
f_{a,b}^{j}(c_{i})=c_{1}$ , $S_{j}=M$ if $c_{1}<f_{a,b}^{j}(c_{i})<c_{2}$, $%
S_{j}=B$ if $f_{a,b}^{j}(c_{i})=c_{2}$ and $S_{j}=R$ if $%
f_{a,b}^{j}(c_{i})>c_{2}$. If we denote by $n_{M}$ the frequency of the
symbol $M$ in a finite subsequence of $S$ we can define the $M$-parity of
this subsequence according to whether $n_{M}$ is even or odd. In what
follows (see \cite{Mi-Th 88}) we define an order relation in $\Sigma
_{5}=\{L,A,M,B,R\}^{\mathbb{N}}$ that depends on the $M$-parity.

Let $V$ be a vector space of three dimension defined over the rationals
having as a basis the formal symbols $\{L,M,R\}$, then to each sequence of
symbols $S=S_{1}S_{2}\dots S_{j}\dots $ we can associate a sequence $\theta
=\theta _{0}\ldots \theta _{j}\ldots $ of vectors from $V$, setting $\theta
_{j}=\prod_{i=0}^{j-1}\epsilon (S_{i})S_{j}$ with $j>0,\theta _{0}=S_{0}$
when $i=0$ and $\epsilon (L)=-\epsilon (M)=\epsilon (R)=1$, where to the
symbols corresponding to the critical points $c_{1}$ and $c_{2}$ we
associate the vector $\frac{L+M}{2}$ and $\frac{M+R}{2}$. Thus $\epsilon
(A)=\epsilon (B)=0$. Choosing then a linear order in the vector space $V$ in
such a way that the base vectors satisfy $L<M<R$ we are able to order the
sequence $\theta $ lexicographically, that is, $\theta <\bar{\theta}$ iff $%
\theta _{0}=\bar{\theta}_{0},\ldots ,\theta _{j-1}=\bar{\theta}_{j-1}$ and $%
\theta _{j}<\bar{\theta}_{j}$ for some integer $i\geq 0$. Finally,
introducing $t$ as an undetermined variable and taking $\theta _{j}$ as the
coefficients of a formal power series $\theta $ (invariant coordinate) we
obtain $\theta =\theta _{0}+\theta _{1}t+\ldots =\sum_{j=0}^{\infty }\theta
_{j}t^{j}$.

The sequences of symbols corresponding to periodic orbits of the critical
points $c_{1}$ and $c_{2}$ are $P=AP_{1}P_{2}\ldots P_{p-1}A\ldots $ and $%
Q=BQ_{1}Q_{2}\ldots Q_{q-1}B\ldots $. In what follows we denote by $%
P^{(p)}=P_{1}P_{2}\ldots P_{p-1}A$ and $Q^{(q)}=Q_{1}Q_{2}\ldots Q_{q-1}B$
the periodic blocks associated to $P$ and $Q$. The realizable itineraries of
the critical points $c_{1}$ and $c_{2}$ for the maps previously defined are
called by kneading sequences \cite{Mi-Th 88}.

\section{Symbolic dynamics for symmetric bimodal maps}

Denote by $\mathcal{F}_{KS}$ the set of pairs of kneading sequences $(P,Q)$,
with $(P,Q)$ either a pair of stable orbits, or a doubly stable orbit. In
Table 1, we give the subset of kneading sequences, with length $p,q<5$.

\begin{center}
Table 1. Kneading data for bimodal maps (detail)
\end{center}

$\hspace*{-0.5cm}{\tiny 
\begin{tabular}{|l|l|l|l|l|l|l|l|l|l|l|l|l|l|l|l|l|l|l|l|l|l|}
\hline\hline
& 1 & 2 & 3 & 4 & 5 & 6 & 7 & 8 & 9 & 10 & 11 & 12 & 13 & 14 & 15 & 16 & 17
& 18 & 19 & 20 & 21 \\ \hline\hline
{1} &  &  &  &  &  &  &  &  &  &  &  &  &  &  &  &  &  &  &  & * &  \\ \hline
{2} &  &  &  &  &  &  &  &  &  &  &  &  &  &  &  & * &  & * &  & * &  \\ 
\hline
{3} &  &  &  &  &  &  &  & * &  & * &  &  &  &  &  & * &  & * &  & * &  \\ 
\hline
{4} &  &  &  &  &  &  & * &  & * &  &  &  & * &  &  &  & * &  & * &  & * \\ 
\hline
{5} &  &  &  &  & $\circledast $ & * &  & * &  & * &  &  &  &  &  & * &  & *
&  & * &  \\ \hline
{6} &  &  &  &  & * & $\circledast $ &  & * &  & * &  &  &  &  &  & * &  & *
&  & * &  \\ \hline
{7} &  &  &  & * &  &  & $\circledast $ &  & * &  & * &  & * &  & * &  & * & 
& * &  & * \\ \hline
{8} &  &  & * &  & * & * &  & $\circledast $ &  & * &  &  &  &  &  & * &  & *
&  & * &  \\ \hline
{9} &  &  &  & * &  &  & * &  & $\circledast $ &  & * &  & * &  & * &  & * & 
& * &  & * \\ \hline
{10} &  &  & * &  & * & * &  & * &  & $\circledast $ &  &  &  &  &  & * &  & 
* &  & * &  \\ \hline
{11} &  &  &  &  &  &  & * &  & * &  & $\circledast $ &  & * &  & * &  & * & 
& * &  & * \\ \hline
{12} &  &  &  &  &  &  &  &  &  &  &  &  &  &  &  & * &  & * &  & * &  \\ 
\hline
{13} &  &  &  & * &  &  & * &  & * &  & * &  & $\circledast $ &  & * &  & *
&  & * &  & * \\ \hline
{14} &  &  &  &  &  &  &  &  &  &  &  &  &  & $\circledast $ &  & * &  & * & 
& * &  \\ \hline
{15} &  &  &  &  &  &  & * &  & * &  & * &  & * &  & $\circledast $ &  & * & 
& * &  & * \\ \hline
{16} &  & * & * &  & * & * &  & * &  & * &  & * &  & * &  & $\circledast $ & 
& * &  & * &  \\ \hline
{17} &  &  &  & * &  &  & * &  & * &  & * &  & * &  & * &  & $\circledast $
&  & * &  & * \\ \hline
{18} &  & * & * &  & * & * &  & * &  & * &  & * &  & * &  & * &  & $%
\circledast $ &  & * &  \\ \hline
{19} &  &  &  & * &  &  & * &  & * &  & * &  & * &  & * &  & * &  & $%
\circledast $ &  & * \\ \hline
{20} & * & * & * &  & * & * &  & * &  & * &  & * &  & * &  & * &  & * &  & $%
\circledast $ &  \\ \hline
{21} &  &  &  & * &  &  & * &  & * &  & * &  & * &  & * &  & * &  & * &  & $%
\circledast $ \\ \hline\hline
\end{tabular}%
}$

\bigskip

Legend: For the lines of the table, we have:

\begin{tabular}{||c||}
\hline\hline
{\tiny 1-RLLA, 2-RLA, 3-RLMA, 4-RLB, 5-RA, 6-RMRA, 7-RMB, 8-RMMA,
9-RMMB,10-RMA,11-RMLB,} \\ \hline
{\tiny 12-RMLA, 13-RB, 14-RRLA, 15-RRLB, 16-RRA, 17-RRMB, 18-RRMA, 19-RRB,
20-RRRA, 21-RRRB} \\ \hline\hline
\end{tabular}

\bigskip

The corresponding columns are given by the conjugate of the previous
sequences.

We define a tree $\mathcal{D}$ that corresponds to the diagonal in $\mathcal{%
F}_{KS}$ and codify the symmetric bimodal maps. Each element $S\in \mathcal{D%
}$ is from one of the following types: $S$ is a pair of stable orbits, i.e., 
$S=(P,\overline{P})=(P^{(p-1)}A,\overline{P^{(p-1)}}B)$; otherwise, $S$ is a
doubly stable orbit, i.e., $S=P^{(p-1)}B\overline{P^{(p-1)}}A$, where $%
\overline{P^{(p-1)}}=\overline{P_{1}}\overline{P_{2}}...\overline{P_{p-1}}$
with $\overline{P_{i}}=R$ if $P_{i}=L$, $\overline{P_{i}}=M$ if $P_{i}=M$
and $\overline{P_{i}}=L$ if $P_{i}=R$, and $1\leq i\leq p-1.$

Note that the set $\mathcal{D}$ is ordered considering the order of the
sequences $P$ (or the inverse order in $\bar{P}$) induced by the order on
the symbols $-R<-B<-M<-A<-L<L<A<M<B<R$. 
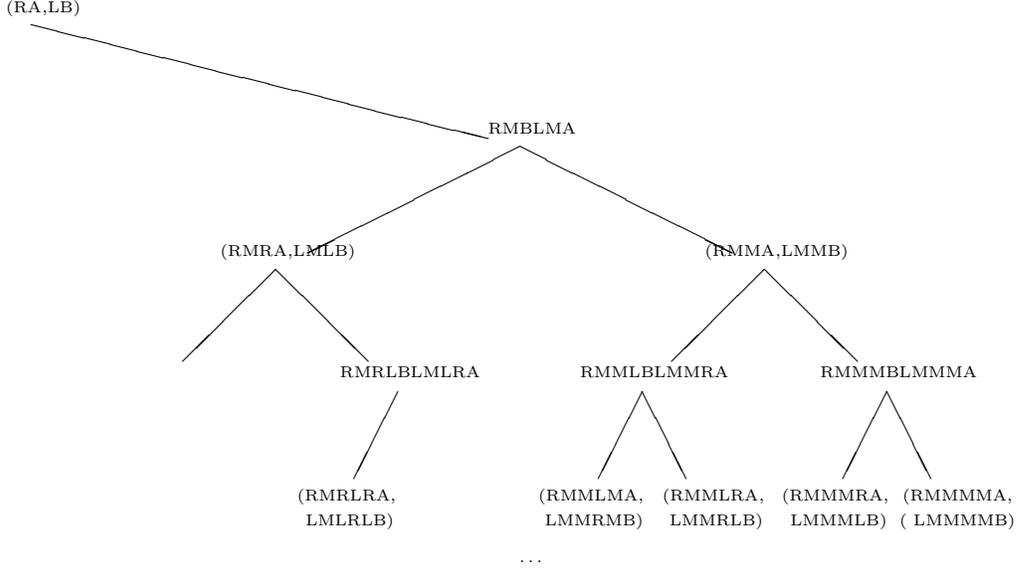
\begin{figure}[th]
\setlength{\unitlength}{0.65mm} {\tiny 
\begin{picture}(180,110)(20,0) 
\put(0,110){\makebox(5,5)[t]{(RA,LB)}}
\put(0,110){\line(4,-1){93.5}} \put(100,85) {\makebox(5,5)[t]{RMBLMA}}
\put(100,85){\line(-2,-1){43.5}} \put(50,60) {\makebox(5,5)[t]{(RMRA,LMLB)}}
\put(100,85){\line(2,-1){43.5}} \put(150,60) {\makebox(5,5)[t]{(RMMA,LMMB)}}
\put(50,60){\line(-1,-1){19}} \put(50,60){\line(1,-1){19}} \put(75,35)
{\makebox(5,5)[t]{RMRLBLMLRA}} \put(150,60){\line(-1,-1){19}} \put(125,35)
{\makebox(5,5)[t]{RMMLBLMMRA}} \put(150,60){\line(1,-1){19}} \put(175,35)
{\makebox(5,5)[t]{RMMMBLMMMA}} \put(75,35){\line(-1,-2){9}} \put(62,10)
{\makebox(5,5)[t]{(RMRLRA,}} \put(62, 5) {\makebox(5,5)[t]{ LMLRLB)}}
\put(125,35){\line(-1,-2){9}} \put(112,10) {\makebox(5,5)[t]{(RMMLMA,}}
\put(112,5) {\makebox(5,5)[t]{ LMMRMB)}} \put(125,35){\line(1,-2){9}}
\put(137,10) {\makebox(5,5)[t]{(RMMLRA,}} \put(137,5) {\makebox(5,5)[t]{
LMMRLB)}} \put(175,35){\line(-1,-2){9}} \put(162,10)
{\makebox(5,5)[t]{(RMMMRA,}} \put(162,5) {\makebox(5,5)[t]{ LMMMLB)}}
\put(175,35){\line(1,-2){9}} \put(187,10) {\makebox(5,5)[t]{(RMMMMA,}}
\put(187,5) {\makebox(5,5)[t]{( LMMMMB)}}
\put(100,0){\makebox(5,5)[b]{\ldots }} 
\end{picture}}
\caption{The Tree $\mathcal{D}_{1}.$}
\label{fig:aci}
\end{figure}

Let $\mathcal{D}_{1}$ a subset of $\mathcal{D}$ with elements between $%
(M^{\infty },M^{\infty })$ and $(RM^{\infty },LM^{\infty })$, see Figure \ref%
{fig:aci}. Let $S^{(2p)}=(P^{(p-1)}A,\overline{P^{(p-1)}}B)$ or $%
S^{(2p)}=P^{(p-1)}B\overline{P^{(p-1)}}A$ and consider a full tree $\mathcal{%
T}$ which its elements are also between $(M^{\infty },M^{\infty })$ and $%
(RM^{\infty },LM^{\infty })$ and characterized by each vertex branch in two
edges following the next rule:

Alternatively the vertices in each level of the tree are doubly stable 
\newline
$P^{(p-1)}B\overline{P^{(p-1)}}A$ or pairs of stable orbits $(P^{(p-1)}A,%
\overline{P^{(p-1)}}B)$. The doubly stable orbits occur in odd levels and
the pairs of stable orbits in even levels. For the doubly stable orbit $%
P^{(p-1)}B\overline{P^{(p-1)}}A$ and according to the $M-$parity of $%
P^{(p-1)}$ is even or odd than the branching order can be described
respectively by: 
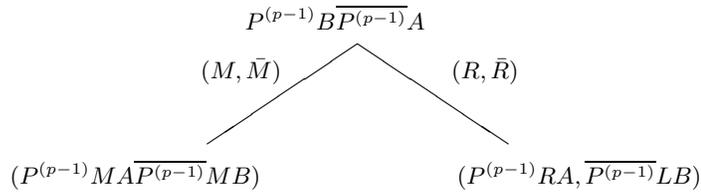
\begin{equation*}
\setlength{\unitlength}{1mm}{\small \begin{picture}(200,30)(0,0)
\put(57,26){\makebox(20,5)[c]{$P^{(p-1)}B\overline{P^{(p-1)}}A$}}
\put(70,25){\line(-3,-2){20}} \put(52,18) {\makebox(5,5)[t]{$(M,\bar M)$}}
\put(31,5){\makebox(20,5)[c]{$(P^{(p-1)}MA\overline{P^{(p-1)} }MB)$ }}
\put(70,25){\line(3,-2){20}} \put(85,18) {\makebox(5,5)[t]{$(R,\bar R)$ }}
\put(90,5) {\makebox(20,5)[c]{$(P^{(p-1)}RA,\overline{P^{(p-1)} }LB) $ }}
\end{picture}}
\end{equation*}%
\begin{equation*}
\setlength{\unitlength}{1mm}{\small 
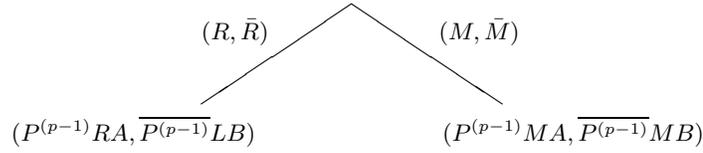
\begin{picture}(200,25)(0,0)
\put(57,26){\makebox(20,5)[c]{$P^{(p-1)}B \overline{P^{(p-1)} }A$ }}
\put(70,25){\line(-3,-2){20}} \put(52,18) {\makebox(5,5)[t]{$(R,\bar R) $}}
\put(31,5){\makebox(20,5)[t]{ $(P^{(p-1)}RA, \overline{P^{(p-1)} }LB) $ }}
\put(70,25){\line(3,-2){20}} \put(85,18) {\makebox(5,5)[t]{$(M,\bar M)$ }}
\put(90,5) {\makebox(20,5)[t]{$(P^{(p-1)}MA, \overline{P^{(p-1)} } MB) $ }}
\end{picture}}
\end{equation*}%
For the pairs of stable orbits the branching order can be described by: 
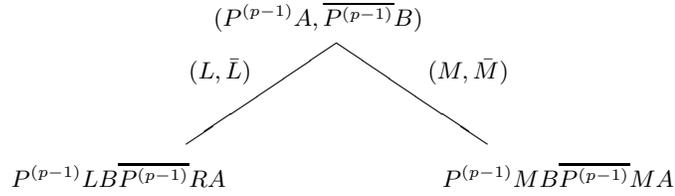
\begin{equation*}
\setlength{\unitlength}{1mm}{\small \begin{picture}(250,35)(0,0)
\put(57,26){\makebox(20,5)[c]{ $(P^{(p-1)}A, \overline{P^{(p-1)} }B)$}}
\put(70,25){\line(-3,-2){20}} \put(52,18) {\makebox(5,5)[t]{ $(L,\bar L)$ }}
\put(31,5){\makebox(20,5)[c]{ $P^{(p-1)}LB \overline{P^{(p-1)} }RA$ }}
\put(70,25){\line(3,-2){20}} \put(85,18) {\makebox(5,5)[t]{ $(M,\bar M)$ }}
\put(90,5) {\makebox(20,5)[c]{$P^{(p-1)}MB \overline{P^{(p-1)} }MA$ }}
\end{picture}}
\end{equation*}%
\begin{equation*}
\setlength{\unitlength}{1mm}{\small 
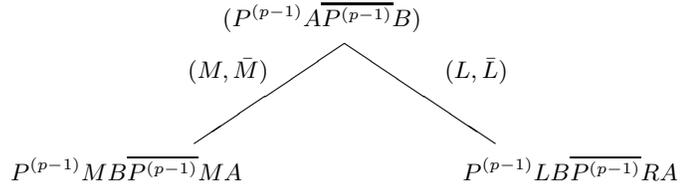
\begin{picture}(250,35)(0,0)
\put(57,26){\makebox(20,5)[c]{$(P^{(p-1)}A \overline{P^{(p-1)} }B)$}}
\put(70,25){\line(-3,-2){20}} \put(52,18){\makebox(5,5)[t]{$(M,\bar M)$}}
\put(31,5){\makebox(20,5)[t]{$P^{(p-1)}MB \overline{P^{(p-1)} }MA$}}
\put(70,25){\line(3,-2){20}} \put(85,18){\makebox(5,5)[t]{$(L,\bar L)$}}
\put(90,5){\makebox(20,5)[t]{$P^{(p-1)}LB \overline{P^{(p-1)} }RA$}}
\end{picture}}
\end{equation*}%
according to the $M-$parity of $P^{(p-1)}A$ is respectively even or odd.

\begin{figure}[th]
\setlength{\unitlength}{0.6mm} {\footnotesize 
\begin{picture}(180,100)(20,0) 
\put(100,75) {\makebox(5,5)[t]{}}
\put(100,75){\line(-2,-1){43.5}} \put(50,50) {\makebox(5,5)[t]{M}}
\put(100,75){\line(2,-1){43.5}} \put(150,50) {\makebox(5,5)[t]{R}}
\put(50,50){\line(-1,-1){19}} \put(25,25) {\makebox(5,5)[t]{MM}}
\put(50,50){\line(1,-1){19}} \put(75,25) {\makebox(5,5)[t]{ML}}
\put(150,50){\line(-1,-1){19}} \put(125,25) {\makebox(5,5)[t]{RL}}
\put(150,50){\line(1,-1){19}} \put(175,25) {\makebox(5,5)[t]{RM}}
\put(75,25){\line(-1,-2){9}} \put(62,0) {\makebox(5,5)[t]{MLR}}
\put(125,25){\line(-1,-2){9}} \put(112,0) {\makebox(5,5)[t]{RLM}}
\put(125,25){\line(1,-2){9}} \put(137,0) {\makebox(5,5)[t]{RLR}}
\put(175,25){\line(-1,-2){9}} \put(162,0) {\makebox(5,5)[t]{RMR}}
\put(175,25){\line(1,-2){9}} \put(187,0) {\makebox(5,5)[t]{RMM}}
\put(25,25){\line(1,-2){9}} \put(8,0) {\makebox(5,5)[t]{MMM}}
\put(25,25){\line(-1,-2){9}} \put(42,0) {\makebox(5,5)[t]{MMR}}
\put(75,25){\line(1,-2){9}} \put(88,0) {\makebox(5,5)[t]{MLM}}
\put(100,-5){\makebox(5,5)[b]{\bf \ldots }} \end{picture}}
\caption{The Tree $\mathcal{T}$}
\label{figfull}
\end{figure}
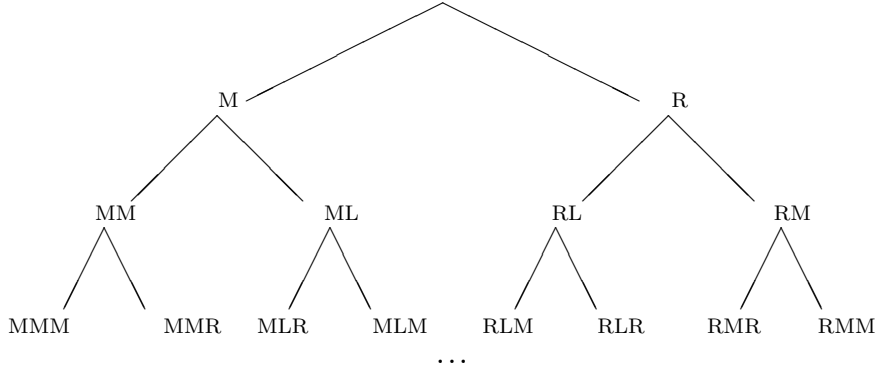

Using these rules we get, as mentioned before, the full tree $\mathcal{T}$,
see Figure \ref{figfull}. 
The next result establish that to each $S=(P^{(p-1)}A,\overline{P^{(p-1)}}B)$%
, or \newline
$P^{(p-1)}B\overline{P^{(p-1)}}A$ in $\mathcal{D}_{1}$, corresponds a
sequence $P^{(p-1)}$ in $\mathcal{T}$.

\begin{lemma}
If $S\in \mathcal{D}_{1}$, then $P^{(p-1)}\in \mathcal{T}$.
\end{lemma}

\begin{proof}
Let $S=P^{(p-1)}B\overline{P^{(p-1)} }A\in \mathcal{D}_{1}$ be a doubly
stable orbit (odd level), with odd $M$-parity. Then, we have: 
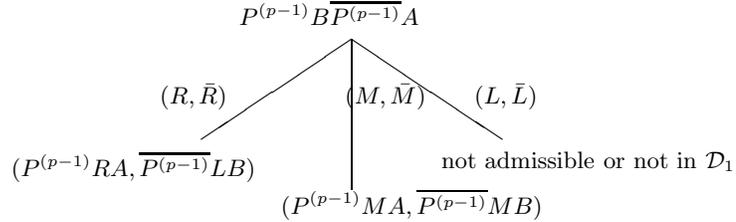
\begin{equation*}
\setlength{\unitlength}{1mm}{\small \begin{picture}(250,30)(0,0)
\put(57,26){\makebox(20,5)[c]{ $ P^{(p-1)} B\overline{P^{(p-1)} }A$ }}
\put(70,25){\line(-3,-2){20}} \put(47,14){\makebox(5,5)[t]{$ (R,\bar R)$ }}
\put(31,5){\makebox(20,5)[t]{ $(P^{(p-1)}RA,\overline{P^{(p-1)} }LB)$ }}
\put(70,25){\line(0, -2){20}} \put(72,14){\makebox(5,5)[t]{ $(M,\bar M)$ }}
\put(68,0){\makebox(20,5)[t]{ $(P^{(p-1)}MA, \overline{P^{(p-1)} }MB) $ }}
\put(70,25){\line(3,-2){20}} \put(88,14){\makebox(5,5)[t]{ $(L,\bar L)$ }}
\put(92,5){\makebox(20,5)[t]{not admissible or not in ${\mathcal D}_1$ }}
\end{picture}}
\end{equation*}%
The doubly stable orbit $P^{(p-1)}B\overline{P^{(p-1)} }A$ leads, on the
next level, to the pairs of stable orbits given by $(P^{(p-1)}XA, \overline{%
P^{(p-1)} } \bar{X}B)=(RM\ldots XA,LM\ldots \bar{X}B)$. Note that when $(X,%
\bar{X})=(L,\bar{L})$ then 
\begin{equation*}
\sigma ^{(p-1)}(LM\ldots \bar{X}B)=\sigma ^{(p-1)}(LM\ldots RB)=RB\ldots
>RM\ldots
\end{equation*}%
which is not admissible or is not in $\mathcal{D}_{1}$. In the same way, the
doubly stable ones obtained from pairs of stable orbits follows the rule in $%
\mathcal{T}$ because now the branch associated to $(R,\bar{R})$ is not
admissible. The proof is analogous for the case when the $M$-parity of $%
P^{(p-1)}A$ is even.
\end{proof}

In what follows we denote by $\mathcal{T}_{KS}$ the set of kneading
sequences associated to unimodal maps. Then, we have:

\begin{theorem}
The tree $\mathcal{D}_{1}$ is isomorphic to $\mathcal{T}_{KS}$.
\end{theorem}

\begin{proof}
Let $\mathcal{E}$ be a complete tree with two symbols $\{L,R\}$ where we
consider the $R$-parity. There exist an isomorphism between $\mathcal{T}$
and $\mathcal{E}$, where each symbol $L$ in $\mathcal{E}$ corresponds to a
symbol $M$ in $\mathcal{T}$ and each symbol $R$ in $\mathcal{E}$ corresponds
to a symbol $L$ or $R$ in $\mathcal{T}$ according to the ($k-1$)- level is
even or odd, respectively. Thus the $R$-parity in $\mathcal{E}$ corresponds
to the ($R+L$)-parity in $\mathcal{T}$ and so, we have:

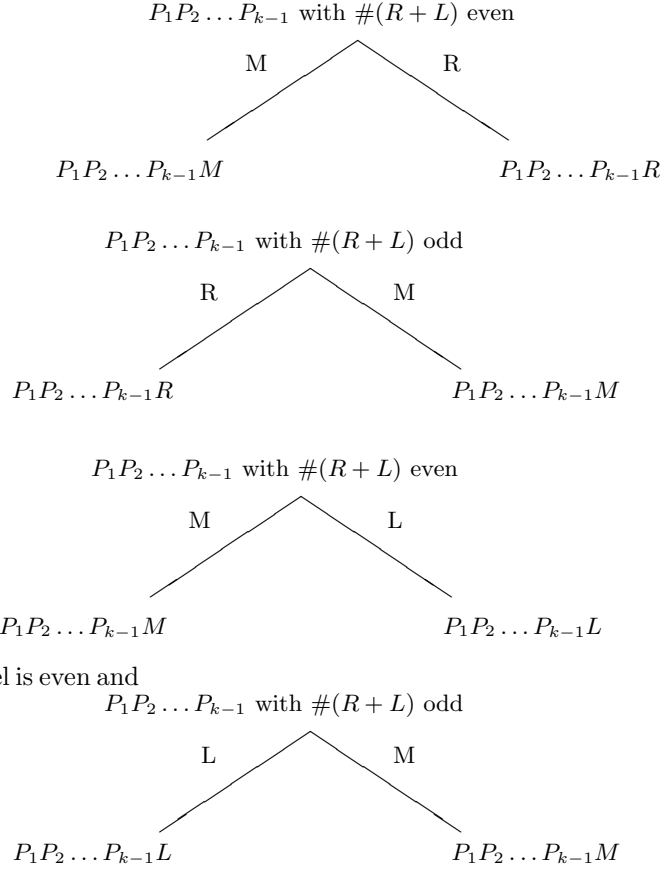
\setlength{\unitlength}{1mm} {\small 
\begin{picture}(200,30)(10,0)
      \put(57,26){\makebox(20,5)[c]{$P_1P_2\ldots P_{k-1}$  
 with $\# (R+L)$ even }}
      \put(70,25){\line(-3,-2){20}}
      \put(54,18){\makebox(5,5)[t]{M}}
      \put(31,5){\makebox(20,5)[c]{$P_1P_2\ldots P_{k-1}M$}}
      \put(70,25){\line(3,-2){20}}
      \put(80,18){\makebox(5,5)[t]{R}}
      \put(90,5){\makebox(20,5)[c]{$P_1P_2\ldots P_{k-1}R$ }}
        \end{picture}} \setlength{\unitlength}{1mm} {\small 
\begin{picture}(200,30)(10,0)
    \put(57,26){\makebox(20,5)[c]{$P_1P_2\ldots P_{k-1}$   with $\# (R+L)$ odd }}
     \put(70,25){\line(-3,-2){20}}
     \put(54,18) {\makebox(5,5)[t]{R}}
      \put(31,5){\makebox(20,5)[t]{$P_1P_2\ldots P_{k-1}R$}}
      \put(70,25){\line(3,-2){20}}
     \put(80,18) {\makebox(5,5)[t]{M}}
      \put(90,5) {\makebox(20,5)[t]{$P_1P_2\ldots P_{k-1}M $}}
 \end{picture}} if ($k-1$)-level is even and \setlength{\unitlength}{1mm} 
{\small 
\begin{picture}(200,30)(50,0)
      \put(57,26){\makebox(20,5)[c]{$P_1P_2\ldots P_{k-1}$   with $\# (R+L)$ even }}
      \put(70,25){\line(-3,-2){20}}
      \put(54,18){\makebox(5,5)[t]{M}}
      \put(31,5){\makebox(20,5)[c]{$P_1P_2\ldots P_{k-1}M$}}
      \put(70,25){\line(3,-2){20}}
      \put(80,18){\makebox(5,5)[t]{L}}
      \put(90,5){\makebox(20,5)[c]{$P_1P_2\ldots P_{k-1}L$ }}
     \end{picture}} \setlength{\unitlength}{1mm} {\small 
\begin{picture}(200,30)(10,0)
    \put(57,26){\makebox(20,5)[c]{$P_1P_2\ldots P_{k-1}$   with $\# (R+L)$ odd }}
     \put(70,25){\line(-3,-2){20}}
     \put(54,18) {\makebox(5,5)[t]{L}}
      \put(31,5){\makebox(20,5)[t]{$P_1P_2\ldots P_{k-1}L$}}
      \put(70,25){\line(3,-2){20}}
     \put(80,18) {\makebox(5,5)[t]{M}}
      \put(90,5) {\makebox(20,5)[t]{$P_1P_2\ldots P_{k-1}M $}}
\end{picture}} if ($k-1$)-level is odd. To each admissible vertex ${P}%
^{(p-1)}C$ when we joint $C$ to the end of a block ${P}^{(p-1)}$ in $%
\mathcal{E}$ corresponds the symbol $A$ (or $B$) in the even or odd level in 
$\mathcal{T}$. Thus, to each admissible vertex ${P}^{(p-1)}C$ in $\mathcal{T}%
_{KS}$ corresponds an admissible vertex (${\tilde{P}}^{(p-1)}A,\widetilde{%
\overline{P}}^{(p-1)}B)$ or ${\tilde{P}}^{(p-1)}B{\bar{P}}^{(p-1)}A$ in $%
\mathcal{D}_{1}$. Note that the admissibility in $\mathcal{E}$ corresponds
to the admissibility in $\mathcal{T}$ since the $R$-parity in $\mathcal{E}$
corresponds to the $(R+L)$-parity in $\mathcal{T}$ and the shift $\sigma $
acting in ${P}^{(p-1)}C$ corresponds in $\mathcal{T}$ to a shift $\sigma $
acting in ${\tilde{P}}^{(p-1)}A$ or ${\tilde{P}}^{(p-1)}B{\bar{P}}^{(p-1)}A$%
. In this way, if ${\ P}^{(p-1)}C$ is admissible, that is, 
\begin{equation*}
\sigma ^{i}({P}^{(p-1)}C)\leq {P}^{(p-1)}C,\mbox{ for all $i$ }
\end{equation*}%
then, we also have that: 
\begin{equation*}
\sigma ^{i}({\tilde{P}}^{(p-1)}A)\leq {\tilde{P}}^{(p-1)}A,\mbox{ for all
$i$ }
\end{equation*}%
or 
\begin{equation*}
\sigma ^{i}({\tilde{P}}^{(p-1)}B{\bar{P}}^{(p-1)}A)\leq {\tilde{P}}^{(p-1)}B{%
\bar{P}}^{(p-1)}A,\mbox{ for all $i$ }.
\end{equation*}
\end{proof}

Consider now the Markov matrix associated to a sequence $S=\tilde{P}^{(p-1)}B%
\overline{P^{(p-1)}}A$ or $(\tilde{P}^{(p-1)}A,\overline{P^{(p-1)}}B)$ and
denote by $d_{P}(t)$ the characteristic polynomial of the Markov matrix $%
A_{P}$ where $P=P^{(p-1)}C\in \mathcal{T}_{KS}$ and corresponds to $\tilde{P}%
={\tilde{P}}^{(p-1)}X$ in $\mathcal{D}_{1}$, where $X=A$ or $B$. Then the
following result holds:

\begin{proposition}
To each $S=\tilde{P}^{(p-1)}B\overline{P^{(p-1)}}A$ or $(\tilde{P}^{(p-1)}A,%
\overline{P^{(p-1)}}B)\in \mathcal{D}_{1}$ there exist a decomposition of
the matrix $A_{S}$ of the type: 
\begin{equation*}
A_{S}=\left[ 
\begin{array}{lll}
1 & W_{1} & W_{2} \\ 
0 & 0 & A_{P} \\ 
0 & A_{P} & 0%
\end{array}%
\right] .
\end{equation*}
\end{proposition}

\begin{proof}
Let $S=\tilde{P}^{(p-1)}B\tilde{\bar{P}}^{(p-1)}A$ or $(\tilde{P}^{(p-1)}A,%
\tilde{\bar{P}}^{(p-1)}B)\in \mathcal{D}_{1}$ then the Markov partition
associated to $S$ is given by $\mathcal{P}=\mathcal{P}_{1}\cup \mathcal{P}%
_{2}\cup \mathcal{P}_{3}$ where $\mathcal{P}_{1}=\{I_{i}\colon 1\leq i\leq
p-1\}$, $\mathcal{P}_{2}=I_{p}$, $\mathcal{P}_{3}=\{I_{i}\colon p+1\leq
i\leq 2p-1\}$ and $\partial I_{i}=z_{i+1}-z_{i}$. When $S=({\tilde{P}}%
^{(p-1)}B{\ \tilde{\bar{P}}}^{(p-1)}A)$ we have: 
\begin{equation*}
z_{i}\in J_{1}=\{x_{2j}\colon 0\leq j<p\}\mbox{ if }I_{i}\in \mathcal{P}_{1}
\end{equation*}%
or 
\begin{equation*}
z_{i}\in J_{2}=\{x_{2j+1}\colon 0\leq j<p\}\mbox{ if }I_{i}\in \mathcal{P}%
_{3}
\end{equation*}%
where $x_{0}$ (resp. $x_{p}$) corresponds to the critical point $c_{1}$
(resp. $c_{2}$). On the other hand, if $S=(\tilde{P}^{(p-1)}A,\tilde{\bar{P}}%
^{(p-1)}B)$ then: 
\begin{equation*}
z_{i}\in J_{3}=\{x_{2j},y_{2j}\colon 0\leq j<\frac{p-2}{2}\}\mbox{ if }%
I_{i}\in \mathcal{P}_{1}
\end{equation*}%
or 
\begin{equation*}
z_{i}\in J_{4}=\{x_{2j+1},y_{2j+1}\colon 0\leq j<\frac{p-2}{2}\}\mbox{ if }%
I_{i}\in \mathcal{P}_{3},
\end{equation*}%
and in both cases $\mathcal{P}_{2}=\{I_{p}\}$, with $\partial
I_{p}=z_{p+1}-z_{p}$, where $z_{p}=\max \{J_{1}(\mbox { or  }J_{3})\}$ and $%
z_{p+1}=\min \{J_{2}(\mbox { or }J_{4})\}$.Note also that if we look for the
structure of $\mathcal{D}_{1}$ we conclude that: If $S=\tilde{P}^{(p-1)}B%
\bar{P}^{(p-1)}A\in \mathcal{D}_{1}$ then $S_{2i}\in \{L,A,M\}$, $%
S_{2i+1}\in \{M,B,R\}$ with $0\leq i<p$. If $S=(\tilde{P}^{(p-1)}A,\bar{P}%
^{(p-1)}B)\in \mathcal{D}_{1}$ then ${\tilde{P}}_{2i}\in \{L,M\}$, ${\tilde{P%
}}_{2i+1}\in \{M,R\}$, ${\bar{P}}_{2i}\in \{M,R\}$ and ${\bar{P}}_{2i+1}\in
\{L,M\}$ for $0\leq i\leq \frac{p-2}{2}$. Thus, the even points establish a
Markov shift and the odd points establish another Markov shift that is
isomorphic to the previous one according to the symmetry of the cubic (where 
$S\in \mathcal{D}_{1}$). So, we only have to prove that each one of these
shifts are isomorphic to the unimodal map associated $(A_{P^{(p)}},\sigma
_{A})$. Note that, the even points are all smaller then the fixed point
(that corresponds to the sequence of symbols $M^{\infty }$) whereas the odd
points are all higher then the fixed point. So, we get two unimodal maps
with critical points $c_{1}$ and $c_{2}$ given by sequences of symbols in $%
\{L,A,M\}$ or $\{M,B,R\}$ and by the admissibility unimodal rules. Thus, the
partitions $\mathcal{P}_{1}$ and $\mathcal{P}_{3}$ are equivalent to the
unimodal map associated and so they have the same Markov shifts. Finally, $%
\mathcal{P}_{2}$ introduce a state that transit to itself and also has
transitions for other states that correspond to the transient part of the
dynamics, $W$.
\end{proof}

\begin{corollary}
To each $S=\tilde{P}^{(p-1)}B\overline{P^{(p-1)}}A$ or $(\tilde{P}^{(p-1)}A,%
\overline{P^{(p-1)}}B)\in \mathcal{D}_{1}$ there exist a decomposition of
the characteristic polynomial, $d_{S}(t)=\det (I-t.A_{S}),$ associated to $%
A_{S}$ that is given by: 
\begin{equation*}
d_{S}(t)=(1-t)d_{P}(t)d_{P}(-t)
\end{equation*}%
where $d_{P}(t)=\det (I-t.A_{P}).$
\end{corollary}

\begin{proof}
Note that the decomposition of the characteristic polynomial follows from
the previous decomposition of the Markov matrix.
\end{proof}

\section{$\star $-Product operator}

For the unimodal case it was defined the $\star $-product operator of
symbolic sequences (see \cite{De-Ge-Po 78}). This product turns out to be a
very useful tool to understand the properties of such maps.

In what follows, we extend the $\star $-product operator for the case of
symbolic sequences associated to symmetric bimodal maps. Note that in $%
\mathcal{D}$ the $\star $-operation is consistent with the initial
definition of the $\star $-product introduced by Derrida-Gervois-Pomeau for
the unimodal case (see also \cite{Ma-Tr 88}, \cite{Ll-Mu 91}, \cite{Pe-Zh 00}
for the bimodal case).

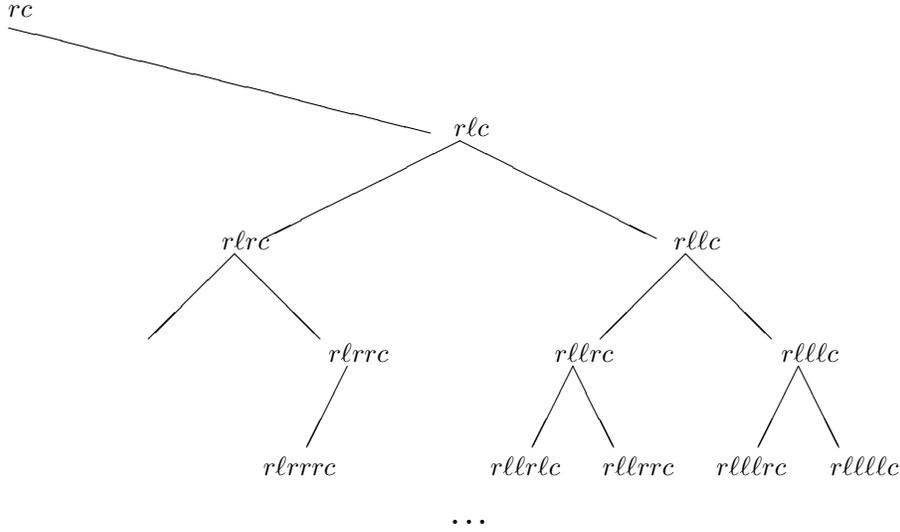
\begin{figure}[th]
\setlength{\unitlength}{0.6mm} {\ 
\begin{picture}(200,120)(20,0) 
\put(0,110){\makebox(5,5)[t]{$rc$}}
\put(0,110){\line(4,-1){93.5}} 
\put(100,85) {\makebox(5,5)[t]{$r\ell c$}}
\put(100,85){\line(-2,-1){43.5}} 
\put(50,60) {\makebox(5,5)[t]{$r\ell rc$}}
\put(100,85){\line(2,-1){43.5}} 
\put(150,60) {\makebox(5,5)[t]{$r\ell \ell c$}}
\put(50,60){\line(-1,-1){19}} 
\put(50,60){\line(1,-1){19}} 
\put(75,35){\makebox(5,5)[t]{$r\ell rrc$}}
 \put(150,60){\line(-1,-1){19}} 
\put(125,35){\makebox(5,5)[t]{$r\ell \ell rc$}} 
\put(150,60){\line(1,-1){19}} 
\put(175,35){\makebox(5,5)[t]{$r\ell \ell \ell c$}} 
\put(75,35){\line(-1,-2){9}} 
\put(62,10){\makebox(5,5)[t]{$r\ell rrrc$}} 
\put(62, 5) {\makebox(5,5)[t]{}}
\put(125,35){\line(-1,-2){9}}
 \put(112,10) {\makebox(5,5)[t]{$r\ell \ell r\ell c$}}
\put(112,5) {\makebox(5,5)[t]{}} 
\put(125,35){\line(1,-2){9}}
\put(137,10) {\makebox(5,5)[t]{$r\ell \ell rrc$}}
 \put(137,5) {\makebox(5,5)[t]{}} 
\put(175,35){\line(-1,-2){9}} 
\put(162,10){\makebox(5,5)[t]{$r\ell \ell \ell rc$}} 
\put(162,5) {\makebox(5,5)[t]{ }}
\put(175,35){\line(1,-2){9}} 
\put(187,10) {\makebox(5,5)[t]{$r\ell \ell \ell \ell c$}}
\put(187,5) {\makebox(5,5)[t]{}}
\put(100,0){\makebox(5,5)[b]{\bf \ldots }} 
\end{picture}}
\caption{The tree $\mathcal{U}$ of unimodal kneading sequences.}
\label{figu}
\end{figure}

According to Theorem 1, the tree $\mathcal{D}_{1}$ is isomorphic to $%
\mathcal{U}$ (ordered set of unimodal kneading sequences, see Fig. \ref{figu}%
).

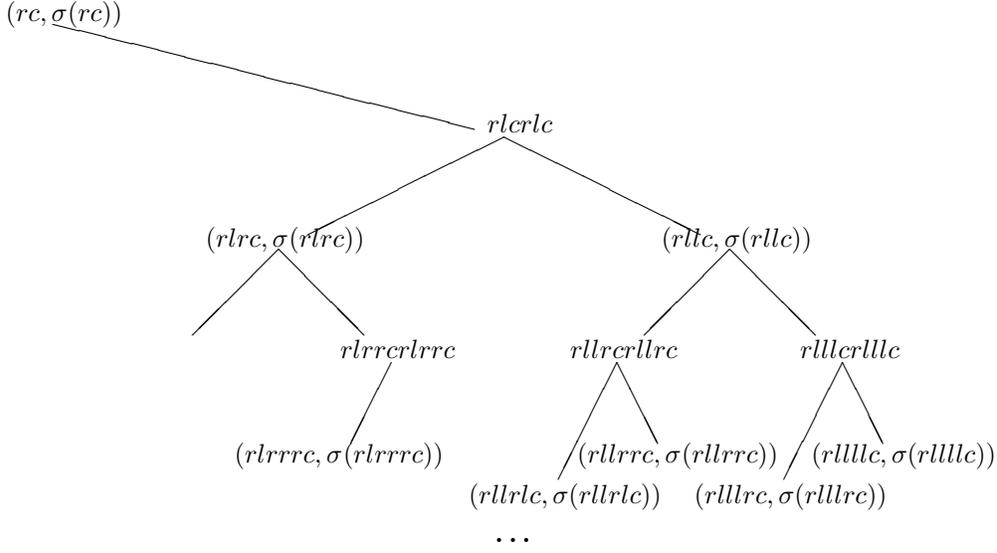
\begin{figure}[th]
\setlength{\unitlength}{0.6mm} {\ 
\begin{picture}(200,120)(20,0) 
\put(0,110){\makebox(5,5)[t]{ $ (rc,\sigma (rc))$ }}
\put(0,110){\line(4,-1){93.5}} \put(100,85) {\makebox(5,5)[t]{ $ rlcrlc $}}
\put(100,85){\line(-2,-1){43.5}} \put(50,60) {\makebox(5,5)[t]{$ (rlrc,\sigma (rlrc))$  }}
\put(100,85){\line(2,-1){43.5}} \put(150,60) {\makebox(5,5)[t]{$ (rllc,\sigma (rllc))$ }}
\put(50,60){\line(-1,-1){19}} 
\put(50,60){\line(1,-1){19}} \put(75,35)
{\makebox(5,5)[t]{$ rlrrcrlrrc$ }}
 \put(150,60){\line(-1,-1){19}} \put(125,35)
{\makebox(5,5)[t]{$ rllrcrllrc$ }} 
\put(150,60){\line(1,-1){19}} \put(175,35)
{\makebox(5,5)[t]{$ rlllcrlllc$ }} 
\put(75,35){\line(-1,-2){9}} \put(62,12)
{\makebox(5,5)[t]{$ (rlrrrc,\sigma (rlrrrc) )$ }} 
\put(125,35){\line(-1,-2){13}}
 \put(112,3) {\makebox(5,5)[t]{$(rllrlc, \sigma ( rllrlc)) $ }}
\put(125,35){\line(1,-2){9}}
\put(137,12) {\makebox(5,5)[t]{$ (rllrrc, \sigma (rllrrc)) $ }}
\put(175,35){\line(-1,-2){13}} 
\put(162,3)
{\makebox(5,5)[t]{$ (rlllrc,\sigma  (rlllrc ))$ }} 
\put(175,35){\line(1,-2){9}} 
\put(187,12) {\makebox(5,5)[t]{$( rllllc,\sigma (rllllc) )$ }}
\put(100,-5){\makebox(5,5)[b]{\bf \ldots }} 
\end{picture}}
\caption{The tree $\mathcal{F}$}
\label{figuu}
\end{figure}

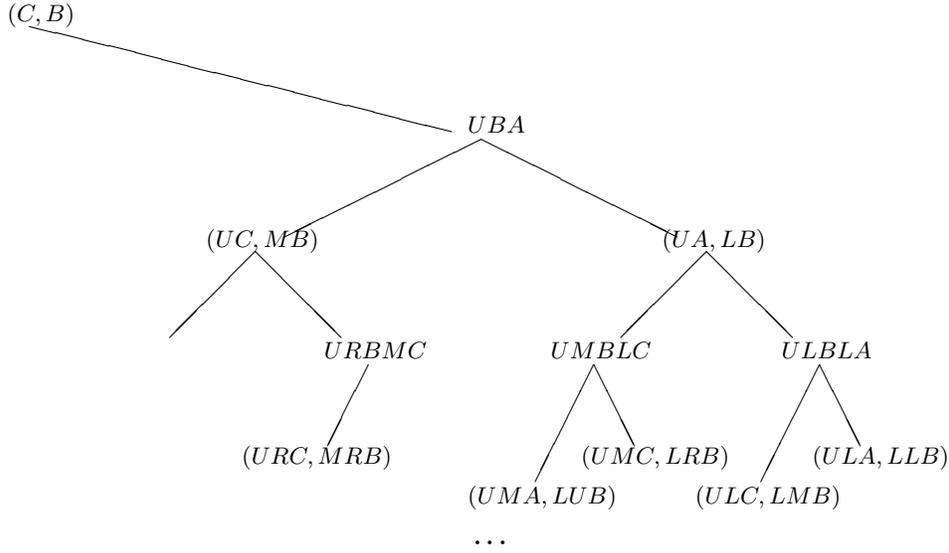
\begin{figure}[th]
\setlength{\unitlength}{0.6mm} {\ {\small 
\begin{picture}(200,120)(20,0) 
\put(0,110){\makebox(5,5)[t]{ $(C, B)$ }}
\put(0,110){\line(4,-1){93.5}} \put(100,85) {\makebox(5,5)[t]{ $UBA$}}
\put(100,85){\line(-2,-1){43.5}} \put(50,60) {\makebox(5,5)[t]{$ (UC, MB)$  }}
\put(100,85){\line(2,-1){43.5}} \put(150,60) {\makebox(5,5)[t]{$ (UA, LB)$ }}
\put(50,60){\line(-1,-1){19}} 
\put(50,60){\line(1,-1){19}} \put(75,35)
{\makebox(5,5)[t]{$ URBMC $ }}
 \put(150,60){\line(-1,-1){19}} \put(125,35)
{\makebox(5,5)[t]{$ UMBLC $ }} 
\put(150,60){\line(1,-1){19}} \put(175,35)
{\makebox(5,5)[t]{$ ULBLA $ }} 
\put(75,35){\line(-1,-2){9}} \put(62,12)
{\makebox(5,5)[t]{$(URC,MRB)$ }} 
\put(125,35){\line(-1,-2){13}}
 \put(112,3) {\makebox(5,5)[t]{$(UMA,LUB) $ }}
\put(125,35){\line(1,-2){9}}
\put(137,12) {\makebox(5,5)[t]{$ (UMC,LRB) $ }}
\put(175,35){\line(-1,-2){13}} 
\put(162,3)
{\makebox(5,5)[t]{$ (ULC,LMB)$ }} 
\put(175,35){\line(1,-2){9}} 
\put(187,12) {\makebox(5,5)[t]{$(ULA,LLB)$ }}
\put(100,-5){\makebox(5,5)[b]{\bf \ldots }} 
\end{picture}
}}
\caption{The tree $\mathcal{G}$ of the second factor for the star product.}
\label{figb}
\end{figure}

From this set we can define another tree $\mathcal{F}=\{\mathcal{F}^{-}=%
\mathcal{UU}$ , if the level is odd, or $\mathcal{F}^{+}=(\mathcal{U},\sigma
(\mathcal{U))}$ if the level is even$\}$, see Fig. \ref{figuu}. 
Now using the symbolic codification applied to $f\circ f$ , with $f$ an
unimodal map, we introduce the following translation rules: 
\begin{equation*}
\ell \ell \longrightarrow L,\ell c\longrightarrow A,\ell r\longrightarrow
M,cr\longrightarrow B,rr\longrightarrow R,rc\longrightarrow C,r\ell
\longrightarrow U.
\end{equation*}%
By applying these rules, Fig. \ref{figuu} can be rewritten as the tree $%
\mathcal{G}=\mathcal{G}^{+}\mathcal{\cup G}^{-}$, see Fig. \ref{figb}. 

\begin{remark}
Let $(x_{1}\dots x_{2n-1}c,\sigma(x_{1}\dots x_{2n-1}c))\in \mathcal{F}^{+}$
or \newline
$x_{1}\dots x_{2n}cx_{1}\dots x_{2n}c \in \mathcal{F}^{-}$, where $x_{i}\in
\{\ell,r\}$, and $\mathcal{F}=\mathcal{F}^{+}\cup \mathcal{F}^{-}$ is the
tree presented in Figure \ref{figuu}. Let $\mathcal{Q}$ be the set of
trimodal kneading data such that the image of both maxima are equal. With $%
D=A$ or $B$, we will write $(P_{1...}P_{2p-1}D,Q_{1}...Q_{2p-1}B,\newline
P_{1}...P_{2p-1}D)\in \mathcal{Q}$, for even levels, and $%
(P_{1...}P_{2p-1}BQ_{1}...Q_{2p-2}A,P_{1}...P_{2p-1}B)\in\mathcal{Q}$, for
odd levels, with $P_{i},Q_{j}\in \{L,A,M,B,R,C,U\}=\{\ell \ell,\ell c,\ell
r,cr,rr,rc,r\ell \}$. With this notation, we can get $\mathcal{G}$ from the
tree $\mathcal{F}$, (see also \cite{La-Se-SR 99b}).
\end{remark}

\noindent Thus, we will consider the following different situations for the
definition of the star product: first, let $F=(P,\overline{P})\in \mathcal{D}
$ and $G=(x_{1} x_{2} \dots x_{x-1}c, x_{1}x_{2} \dots x_{x-1}c)$, with $%
x_{1}x_{2}...x_{x-1} c \in \mathcal{U}$; then, we let $F=PB\overline{P}A \in 
\mathcal{D}$ and $G \in \mathcal{G}$.

\noindent \textbf{Type 1.} Let $F=(P,\overline{P})=(P^{(p-1)}A,\overline{%
P^{(p-1)}}B)\in \mathcal{D}$ be a bimodal kneading data and $G=(X,X)$ be a
pair of unimodal kneading sequences. Then, we have 
\begin{equation*}
F\star G=(P,\overline{P}) \star (X,X)=(P^{(p-1)}\star X^{(x-1)}c,\overline{%
P^{(p-1)}}\star X^{(x-1)}c),
\end{equation*}%
with 
\begin{equation*}
P^{(p-1)}\star X^{(x-1)}c=P^{(p-1)}A_{1}^{\pm }P^{(p-1)}A_{2}^{\pm }\ldots
P_{x-1}^{(p-1)}A_{x-1}^{\pm }P^{(p-1)}A,
\end{equation*}
where 
\begin{equation*}
\begin{array}{ll}
A_{i}^{\pm }=\left\{ 
\begin{array}{l}
M\mbox{  if  }x_{i}=r \\ 
A\mbox{  if  }x_{i}=c \\ 
L\mbox{  if  }x_{i}=\ell%
\end{array}
\right. \text{ if }P\,\,\text{is even, } & A_{i}^{\pm }=\left\{ 
\begin{array}{l}
L\mbox{  if  }x_{i}=r \\ 
A\mbox{  if  }x_{i}=c \\ 
M\mbox{  if  }x_{i}=\ell%
\end{array}
\right. \text{ if }P\,\,\text{is odd.}%
\end{array}%
\end{equation*}
In a similar way, 
\begin{equation*}
\overline{P^{(p-1)}}\star X^{(x-1)}c=\overline{P^{(p-1)}}B_{1}^{\pm }%
\overline{P^{(p-1)}}B_{2}^{\pm }\ldots \overline{P^{(p-1)}}B_{y-1}^{\pm }%
\overline{P^{(p-1)}}B,
\end{equation*}
where 
\begin{equation*}
\begin{array}{ll}
B_{i}^{\pm }=\left\{ 
\begin{array}{l}
M\mbox{  if  }x_{i}=r \\ 
B\mbox{  if  }x_{i}=c \\ 
R\mbox{  if  }x_{i}=\ell%
\end{array}%
\right. \text{ if }\overline{P}\,\,\text{is even, } & B_{i}^{\pm }=\left\{ 
\begin{array}{l}
R\mbox{  if  }x_{i}=r \\ 
B\mbox{  if  }x_{i}=c \\ 
M\mbox{  if  }x_{i}=\ell%
\end{array}
\right. \text{ if }\overline{P}\,\,\text{is odd.}%
\end{array}%
\end{equation*}

\noindent \textbf{Type 2.} Let $F=PB\overline{P}A=P^{(p-1)}B\overline{%
P^{(p-1)}}A\in \mathcal{D}$\ and $G=X^{(x-1)}B{Y}^{(y-1)}D\in \mathcal{G}^{-}
$ (where $D=A$ or $C$) be two kneading data. Then, 
\begin{equation*}
\begin{array}{l}
F\star G = P^{(p-1)}B_{1}^{\pm }\overline{P^{(p-1)}}A_{1}^{%
\pm}P^{(p-1)}B_{2}^{\pm }\overline{P^{(p-1)}}A_{2}^{\pm } \ldots P^{(p-1)}B 
\overline{P^{(p-1)}}A_{x}^{\pm }P^{(p-1)} \\[0.2cm] 
\hspace{50pt} B_{x+1}^{\pm }\overline{P^{(p-1)}}A_{x+1}^{\pm
}P^{(p-1)}B_{x+2}^{\pm }\ldots \overline{P^{(p-1)}}A_{x+y-1}^{\pm
}P^{(p-1)}B_{x+y}^{\pm }\overline{P^{(p-1)}}A.%
\end{array}
\end{equation*}
Let $Z_{i}=$ $X_{i}$, $Z_{x+i}=Y_{i}$ and $Z_{x+y}=D$ then 
\begin{equation*}
P^{(p-1)}B_{i}^{\pm }\overline{P^{(p-1)}}A_{i}^{\pm }=\left\{ 
\begin{array}{l}
P^{(p-1)}M\overline{P^{(p-1)}}M\mbox{  if  }Z_{i}=L \\ 
P^{(p-1)}M\overline{P^{(p-1)}}A\mbox{  if  }Z_{i}=A \\ 
P^{(p-1)}M\overline{P^{(p-1)}}L\mbox{  if  }Z_{i}=M \\ 
P^{(p-1)}B\overline{P^{(p-1)}}L\mbox{  if  }Z_{i}=B \\ 
P^{(p-1)}R\overline{P^{(p-1)}}L\mbox{  if  }Z_{i}=R \\ 
P^{(p-1)}R\overline{P^{(p-1)}}A\mbox{  if  }Z_{i}=C \\ 
P^{(p-1)}R\overline{P^{(p-1)}}M\mbox{  if  }Z_{i}=U%
\end{array}%
\right. \text{ if }P\ \text{and }\overline{P}\,\ \text{are even, }
\end{equation*}%
\begin{equation*}
P^{(p-1)}B_{i}^{\pm }\overline{P^{(p-1)}}A_{i}^{\pm }=\left\{ 
\begin{array}{l}
P^{(p-1)}R\overline{P^{(p-1)}}L\mbox{  if  }Z_{i}=L \\ 
P^{(p-1)}R\overline{P^{(p-1)}}A\mbox{  if  }Z_{i}=A \\ 
P^{(p-1)}R\overline{P^{(p-1)}}M\mbox{  if  }Z_{i}=M \\ 
P^{(p-1)}B\overline{P^{(p-1)}}M\mbox{  if  }Z_{i}=B \\ 
P^{(p-1)}M\overline{P^{(p-1)}}M\mbox{  if  }Z_{i}=R \\ 
P^{(p-1)}M\overline{P^{(p-1)}}A\mbox{  if  }Z_{i}=C \\ 
P^{(p-1)}M\overline{P^{(p-1)}}L\mbox{  if  }Z_{i}=U%
\end{array}%
\right. \text{ if }P\ \text{and }\overline{P}\ \text{are odd, }
\end{equation*}

\noindent \textbf{Type 3.} Let $F=P^{(p-1)}B\overline{P^{(p-1)}}A \in 
\mathcal{D}$ and $G=(X^{(n-1)}D,{Y}^{(n-1)}B)\in \mathcal{G}^{+}$ be two
kneading data. Then 
\begin{eqnarray*}
&& F \star G = (P^{(p-1)}B_{1}^{\pm }\overline{P^{(p-1)}}A_{1}^{%
\pm}P^{(p-1)}B_{2}^{\pm }\overline{P^{(p-1)}}A_{2}^{\pm }\ldots
P^{(p-1)}B_{n}^{\pm }\overline{P^{(p-1)}}A, \\
&& \hspace{20pt} \overline{P^{(p-1)}}A_{n+1}^{\pm }P^{(p-1)}B_{n+1}^{\pm }%
\overline{P^{(p-1)}}A_{n+2}^{\pm }P^{(p-1)}B_{n+2}^{\pm }\ldots \overline{%
P^{(p-1)}}A_{2n}^{\pm }P^{(p-1)}B).
\end{eqnarray*}%
The transformation rules are the same as above for the first sequence 
\begin{equation*}
F\star X^{(x-1)}A=P^{(p-1)}B_{1}^{\pm }\overline{P^{(p-1)}}A_{1}^{\pm
}P^{(p-1)}B_{2}^{\pm }\overline{P^{(p-1)}}A_{2}^{\pm }\ldots
P^{(p-1)}B_{n}^{\pm }\overline{P^{(p-1)}}A,
\end{equation*}%
except that $Z_{i}=B$ cannot occur. For the second position of the pair we
have 
\begin{equation*}
\begin{array}{l}
\sigma ^{p-1}(F)\star \sigma ^{n-1}({Y}^{(n-1)}B)= \\[0.2cm] 
\hspace{30pt} B\overline{P^{(p-1)}}A_{n+1}^{\pm }P^{(p-1)}B_{n+1}^{\pm }%
\overline{P^{(p-1)} }A_{n+2}^{\pm }P^{(p-1)}B_{n+2}^{\pm }\ldots \overline{%
P^{(p-1)}}A_{2n}^{\pm }P^{(p-1)} \\ 
\\ 
\sigma (\sigma ^{p-1}(F)\star \sigma ^{n-1}({Y}^{(n-1)}B))= \\[0.2cm] 
\hspace{30pt} \overline{P^{(p-1)}}A_{n+1}^{\pm }P^{(p-1)}B_{n+1}^{\pm }%
\overline{P^{(p-1)}} A_{n+2}^{\pm }P^{(p-1)}B_{n+2}^{\pm }\ldots \overline{%
P^{(p-1)}}A_{2n}^{\pm}P^{(p-1)}B.%
\end{array}%
\end{equation*}%
where 
\begin{equation*}
B_{i}^{\pm }\overline{P^{(p-1)}}A_{i}^{\pm }P^{(p-1)}=\left\{ 
\begin{array}{l}
M\overline{P^{(p-1)}}MP^{(p-1)}\mbox{  if  }Y_{i}=L \\ 
M\overline{P^{(p-1)}}LP^{(p-1)}\mbox{  if  }Y_{i}=M \\ 
B\overline{P^{(p-1)}}LP^{(p-1)}\mbox{  if  }Y_{i}=B \\ 
R\overline{P^{(p-1)}}LP^{(p-1)}\mbox{  if  }Y_{i}=R \\ 
R\overline{P^{(p-1)}}MP^{(p-1)}\mbox{  if  }Y_{i}=U%
\end{array}%
\right. \text{ if }P\ \text{and }\overline{P}\,\ \text{are even, }
\end{equation*}%
\begin{equation*}
B_{i}^{\pm }\overline{P^{(p-1)}}A_{i}^{\pm }P^{(p-1)}=\left\{ 
\begin{array}{l}
R\overline{P^{(p-1)}}LP^{(p-1)}\mbox{  if  }Y_{i}=L \\ 
R\overline{P^{(p-1)}}MP^{(p-1)}\mbox{  if  }Y_{i}=M \\ 
B\overline{P^{(p-1)}}MP^{(p-1)}\mbox{  if  }Y_{i}=B \\ 
M\overline{P^{(p-1)}}MP^{(p-1)}\mbox{  if  }Y_{i}=R \\ 
M\overline{P^{(p-1)}}LP^{(p-1)}\mbox{  if  }Y_{i}=U%
\end{array}%
\right. \text{ if }P\ \text{and }\overline{P}\,\ \text{are odd, }
\end{equation*}%
The following examples illustrate the definitions given above:

\begin{example}
\begin{equation*}
\mathit{(RMMA,LMMB)\star (r\ell c,r\ell c)=(RMMMRMMLRMMA,LMMMLMMRLMMB).} 
\end{equation*}
\end{example}

\begin{example}
\begin{equation*}
\mathit{RBLA\star ULBLA=RRLMRMLMRBLLRMLMRMLA}. 
\end{equation*}
\end{example}

\begin{example}
\begin{equation*}
\mathit{RBLA\star (UA,LB)=(RRLMRMLA,LLRMLMRB)}. 
\end{equation*}
\end{example}

\begin{remark}
Notice that, regarding Example 3, where both $P$ and $\overline{P}$ are even
sequences, we have, for the second position of the pair, 
\begin{equation*}
\mathit{\sigma(\sigma(RBLA)\star \sigma(LB))=\sigma (BLAR\star BL)= \sigma
(BLLRMLMR)=LLRMLMRB}.
\end{equation*}
\end{remark}

\begin{remark}
The star product in $\mathcal{D}$ is not a true binary product $A \ast B=C,$
for all $A$ and B $\in \mathcal{D}$. When $C$ is factorizable, with $C\in 
\mathcal{D}$ then $B$ must be in $\mathcal{G}$.
\end{remark}

\begin{remark}
Note that for all $A\in \mathcal{D}$ and B $\in \mathcal{G}$. then the
result of the star product defined previous is also in $\mathcal{D}$.
\end{remark}

\section{$\otimes $-Product between Markov matrices}

In the same way, we can extend the $\otimes $-product between Markov
matrices (introduced for unimodal maps in \cite{La-RS-SR 88}) associated to
symmetric bimodal maps, that is 
\begin{equation*}
A_{S}=A_{V}\otimes A_{W}
\end{equation*}%
where $S=V\star W$ with $V\in \mathcal{D}$ and $W\in \mathcal{G}$.

\begin{theorem}
Let $V\in \mathcal{D}{}$ and $W\in \mathcal{G}$ ~then there exists a matrix
product such that 
\begin{equation*}
A_{S}=A_{V\star W}=A_{V}\otimes A_{W}
\end{equation*}
\end{theorem}

\begin{proof}
It is based on a construction of a product on the matrices induced by the $%
\star $-product between kneading sequences. We will give this construction
but only for the $\star $-product between kneading sequences of the first
type. For the others it is technically similar and can be reproduced from
this one. Let $W=(x_{1}x_{2}...x_{k-1}c,x_{1}x_{2}...x_{k-1}c)\in \mathcal{G}%
.$ First of all note that the matrix $A_{W}$ is symmetric and so it can be
written in the form: 
\begin{equation*}
A_{W}=\left[ 
\begin{array}{ccc}
0 & 0 & \widehat{B}_{X} \\ 
N_{1} & 1 & N_{2} \\ 
B_{X} & 0 & 0%
\end{array}%
\right]
\end{equation*}%
where $[N_{1}\quad 1\quad N_{2}]$ is the $k$ row and $[0\quad 1\quad 0]$ is
the $k$ column. Denoting $B_{X}=[b_{ij}]$ then we define the $(k-1)\times
(k-1)$ matrix $\hat{B}_{X}$, by $\hat{B}_{X}=B_{X}$ if $P^{(p)}$ is even and
by $\hat{B}_{X}=[\hat{b}_{mj}]$ with $\hat{b}_{mj}=b_{ij}$ where $m=k-i$ if $%
P^{(p)}$ is odd. Given $V=(P,\overline{P})\in \mathcal{D}$ and $W=(X,%
\overline{X})\in \mathcal{G}$, it is immediate that its associated
transition matrices, $A_{V}$ and $A_{W}$, are square, $(2p-1)$ and $(2k-1)$%
-dimensional matrices, respectively, where $p$ and $k$ denotes the number of
symbols of the sequences $P$ and $X$. Analogously, it is fairly simple to
see that the transition matrix associated with the sequence $V\star W$ is a
square matrix with dimension $(2pk-1)$. Now, we need to show that the
elements of $A_{V\star W}$ are completely determined by the knowledge of the
matrices $A_{V}$ and $A_{W}$. Consider the symbolic shifts of the sequences $%
P\star X$ and $\overline{P}\star \overline{X}$ and denote the corresponding
points of the interval by $p_{i}^{j}$ and $q_{i}^{j}$, that is, $p_{i}^{j}$
will be the point corresponding to the sequence $\sigma
^{p(j-1)+(i-1)}(P\star X)$ and $q_{i}^{j}$ the point corresponding to the
sequence $\sigma ^{p(j-1)+(i-1)}(\overline{P}\star \overline{X})$. When one
considers the collection of points of the interval from all the shifts cited
above, we can see that they appear as groups of blocks of $x$ points.
Considering the order of the shifted sequences $\sigma ^{i}(P)$, $\sigma
^{j}(\overline{P})$, $\sigma ^{k}(X)$ and $\sigma ^{n}(\overline{X})$ and
the way those sequences appear as subsets of the partition induced by the
sequence $V\star W$, we can conclude (see also \cite{La-RS-SR 88}) that the
matrix $A_{V\star W}$ has the following block structure: 
\begin{equation*}
\left[ 
\begin{array}{ccccc}
A_{1,1} & A_{1,2} & \cdots & A_{1,l+m+r} & N_{1} \\ 
\vdots & \vdots &  & \vdots & \vdots \\ 
A_{l,1} & A_{l,2} & \cdots & A_{l,l+m+r} & N_{l} \\ 
Z_{l+1,1} & Z_{l+1,2} & \cdots & Z_{l+1,l+m+r} & \bar{A}_{X} \\ 
A_{l+2,1} & A_{l+2,2} & \cdots & A_{l+2,l+m+r} & N_{l+2} \\ 
\vdots & \vdots &  & \vdots & \vdots \\ 
A_{l+m+1,1} & A_{l+m+1,2} & \cdots & A_{l+m+1,l+m+r} & N_{l+m+1} \\ 
\widetilde{A}_{X} & Z_{l+m+2,2} & \cdots & Z_{l+m+2,l+m+r} & N_{l+m+2} \\ 
N_{l+m+3} & A_{l+m+3,2} & \cdots & A_{l+m+3,l+m+r} & A_{l+m+3,l+m+r+1} \\ 
\vdots & \vdots &  & \vdots & \vdots \\ 
N_{l+m+r+2} & A_{l+m+r+2,2} & \cdots & A_{l+m+r+2,l+m+r} & 
A_{l+m+r+2,l+m+r+1}%
\end{array}
\right],
\end{equation*}
where $l,m$ and $r$ are, respectively, the number of symbols $L,M$ and $R$
in the sequence $V$, and $A_{i,j}$, with $(i,j)\neq (l+m+1,1)$, is either
one of these $k\times k$ matrices 
\begin{eqnarray*}
&&\left[ 
\begin{array}{cccc}
1 & 1 & \cdots & 1 \\ 
0 & 0 & \cdots & 0 \\ 
\vdots & \vdots &  & \vdots \\ 
0 & 0 & \cdots & 0%
\end{array}
\right], \left[ 
\begin{array}{cccc}
1 & 0 & \cdots & 0 \\ 
0 & 1 & \cdots & 0 \\ 
\vdots & \vdots &  & \vdots \\ 
0 & 0 & \cdots & 1%
\end{array}
\right], \\
&& \left[ 
\begin{array}{cccc}
0 & \cdots & 0 & 1 \\ 
0 & \cdots & 1 & 0 \\ 
\vdots &  & \vdots & \vdots \\ 
1 & \cdots & 0 & 0%
\end{array}
\right], \left[ 
\begin{array}{cccc}
0 & 0 & \cdots & 0 \\ 
0 & 0 & \cdots & 0 \\ 
\vdots & \vdots &  & \vdots \\ 
1 & 1 & \cdots & 1%
\end{array}
\right],
\end{eqnarray*}
or a null block, and 
\begin{equation*}
A_{l+m+1,1}= \left[ 
\begin{array}{cccc}
n_{1} & \cdots & n_{y} & n_{y+1} \\ 
b_{1,1} & \cdots & b_{1,y} & m_{1} \\ 
\vdots &  & \vdots & \vdots \\ 
b_{y,1} & \cdots & b_{y,y} & m_{y}%
\end{array}
\right],
\end{equation*}
where the submatrix $[b_{i,j}]$ is defined by the matrix $B_{X}$. The
matrices $Z_{k}$ and $N_{j}$ are null matrices, except $N_{l}$, $N_{l+1}$
and $N_{l+m+2}$ that can contain some elements $1$. The distribution of the
previous blocks $A_{i,j}$, with $(i,j)\neq (l+m+1,1)$, is given by the
structure of the matrix $A_{V}$. On the other hand, the internal structure
of each block $A_{i,j}$ is determined by the order of the shifts of the
sequence $W$. For the case of the block $A_{l+m+1,1}$, its submatrix $%
[b_{i,j}]$ has an internal structure determined by $B_{X}$. The elements $%
n_{i}$ and $m_{j}$ are null except those needed to preserve the continuity
of the transitions (from the fact that $f$ is a continuous function).
Analogously, the block $\bar{A}_{X}$ is determined by $\widehat{B}_{X}$.
Finally, the blocks $N_{l}$, $N_{l+1}$, $N_{l+m+2}$ are null except, once
again, for the elements needed to preserve the continuity of the transitions.
\end{proof}

The following example illustrate the use of the previous theorem.

\begin{example}
Let $\mathit{V=(RMMA,LMMB)}$ and $W=(r\ell c,r\ell c)$. Then, we have $%
\mathit{S=V \star W=(RMMMRMMLRMMA,LMMMLMMRLMMB)}$, and the $\otimes$-product
of its matrices, 
\begin{equation*}
A_{V}=\left[ 
\begin{array}{ccccccc}
0 & 0 & 0 & 0 & 1 & 1 & 1 \\ 
0 & 0 & 0 & 0 & 0 & 0 & 1 \\ 
0 & 0 & 0 & 0 & 0 & 1 & 0 \\ 
0 & 0 & 1 & 1 & 1 & 0 & 0 \\ 
0 & 1 & 0 & 0 & 0 & 0 & 0 \\ 
1 & 0 & 0 & 0 & 0 & 0 & 0 \\ 
1 & 1 & 1 & 0 & 0 & 0 & 0%
\end{array}
\right] ; \qquad A_{W}=\bigl(\left[ 
\begin{array}{cc}
0 & 1 \\ 
1 & 1%
\end{array}
\right] ,\left[ 
\begin{array}{cc}
0 & 1 \\ 
1 & 1%
\end{array}
\right] \bigr), 
\end{equation*}
$A_{S}=A_{V \star W}$, is given by 
\begin{equation*}
\left[ 
\begin{array}{ccccccccccccccccccccccc}
0 & 0 & 0 & 0 & 0 & 0 & 0 & 0 & 0 & 0 & 0 & 0 & 1 & 0 & 0 & 0 & 0 & 0 & 0 & 0
& 0 & 0 & 0 \\ 
0 & 0 & 0 & 0 & 0 & 0 & 0 & 0 & 0 & 0 & 0 & 0 & 0 & 1 & 0 & 0 & 0 & 0 & 0 & 0
& 0 & 0 & 0 \\ 
0 & 0 & 0 & 0 & 0 & 0 & 0 & 0 & 0 & 0 & 0 & 0 & 0 & 0 & 1 & 1 & 1 & 1 & 1 & 1
& 1 & 1 & 0 \\ 
0 & 0 & 0 & 0 & 0 & 0 & 0 & 0 & 0 & 0 & 0 & 0 & 0 & 0 & 0 & 0 & 0 & 0 & 0 & 0
& 0 & 0 & 1 \\ 
0 & 0 & 0 & 0 & 0 & 0 & 0 & 0 & 0 & 0 & 0 & 0 & 0 & 0 & 0 & 0 & 0 & 0 & 0 & 0
& 0 & 1 & 1 \\ 
0 & 0 & 0 & 0 & 0 & 0 & 0 & 0 & 0 & 0 & 0 & 0 & 0 & 0 & 0 & 0 & 0 & 0 & 0 & 0
& 1 & 0 & 0 \\ 
0 & 0 & 0 & 0 & 0 & 0 & 0 & 0 & 0 & 0 & 0 & 0 & 0 & 0 & 0 & 0 & 0 & 0 & 0 & 1
& 0 & 0 & 0 \\ 
0 & 0 & 0 & 0 & 0 & 0 & 0 & 0 & 0 & 0 & 0 & 0 & 0 & 0 & 0 & 0 & 0 & 0 & 1 & 0
& 0 & 0 & 0 \\ 
0 & 0 & 0 & 0 & 0 & 0 & 0 & 0 & 0 & 0 & 0 & 0 & 0 & 0 & 0 & 0 & 0 & 1 & 0 & 0
& 0 & 0 & 0 \\ 
0 & 0 & 0 & 0 & 0 & 0 & 0 & 0 & 0 & 0 & 0 & 0 & 0 & 0 & 0 & 0 & 1 & 0 & 0 & 0
& 0 & 0 & 0 \\ 
0 & 0 & 0 & 0 & 0 & 0 & 0 & 0 & 0 & 0 & 0 & 0 & 0 & 0 & 0 & 1 & 0 & 0 & 0 & 0
& 0 & 0 & 0 \\ 
0 & 0 & 0 & 0 & 0 & 0 & 0 & 0 & 1 & 1 & 1 & 1 & 1 & 1 & 1 & 0 & 0 & 0 & 0 & 0
& 0 & 0 & 0 \\ 
0 & 0 & 0 & 0 & 0 & 0 & 0 & 1 & 0 & 0 & 0 & 0 & 0 & 0 & 0 & 0 & 0 & 0 & 0 & 0
& 0 & 0 & 0 \\ 
0 & 0 & 0 & 0 & 0 & 0 & 1 & 0 & 0 & 0 & 0 & 0 & 0 & 0 & 0 & 0 & 0 & 0 & 0 & 0
& 0 & 0 & 0 \\ 
0 & 0 & 0 & 0 & 0 & 1 & 0 & 0 & 0 & 0 & 0 & 0 & 0 & 0 & 0 & 0 & 0 & 0 & 0 & 0
& 0 & 0 & 0 \\ 
0 & 0 & 0 & 0 & 1 & 0 & 0 & 0 & 0 & 0 & 0 & 0 & 0 & 0 & 0 & 0 & 0 & 0 & 0 & 0
& 0 & 0 & 0 \\ 
0 & 0 & 0 & 1 & 0 & 0 & 0 & 0 & 0 & 0 & 0 & 0 & 0 & 0 & 0 & 0 & 0 & 0 & 0 & 0
& 0 & 0 & 0 \\ 
0 & 0 & 1 & 0 & 0 & 0 & 0 & 0 & 0 & 0 & 0 & 0 & 0 & 0 & 0 & 0 & 0 & 0 & 0 & 0
& 0 & 0 & 0 \\ 
1 & 1 & 0 & 0 & 0 & 0 & 0 & 0 & 0 & 0 & 0 & 0 & 0 & 0 & 0 & 0 & 0 & 0 & 0 & 0
& 0 & 0 & 0 \\ 
1 & 0 & 0 & 0 & 0 & 0 & 0 & 0 & 0 & 0 & 0 & 0 & 0 & 0 & 0 & 0 & 0 & 0 & 0 & 0
& 0 & 0 & 0 \\ 
0 & 1 & 1 & 1 & 1 & 1 & 1 & 1 & 1 & 0 & 0 & 0 & 0 & 0 & 0 & 0 & 0 & 0 & 0 & 0
& 0 & 0 & 0 \\ 
0 & 0 & 0 & 0 & 0 & 0 & 0 & 0 & 0 & 1 & 0 & 0 & 0 & 0 & 0 & 0 & 0 & 0 & 0 & 0
& 0 & 0 & 0 \\ 
0 & 0 & 0 & 0 & 0 & 0 & 0 & 0 & 0 & 0 & 1 & 0 & 0 & 0 & 0 & 0 & 0 & 0 & 0 & 0
& 0 & 0 & 0%
\end{array}
\right]
\end{equation*}
\end{example}


\end{document}